\documentstyle{article}
\begin{document}

\begin{center}
{\LARGE
{\bf Gradient like Morse-Smale

\medskip
dynamical systems on 4-manifolds}}

\bigskip
{\Large
{\bf A.O.Prishlyak}
}
\end{center}

{\large
\bigskip
{\it Kiev University, Geometry Department,

Ukraine, 252083, Kiev, Xytopcka 3}

{\it e-mail:} prish@mechmat.univ.kiev.ua

\medskip
\noindent
UDK 517.91

\noindent
AMS MSC 1991: Primary 58F09. Secondary 57N13, 58F25.}

\bigskip

Abstract. The complete invariant for gradient like Morse-Smale
dynamical systems (vector fields and diffeomorphisms) on closed
4-manifolds are constructed. It is same
as Kirby diagram in a case of polar vector field without fixed points of
index 3.

{\large
\bigskip
\noindent
{\bf 1. Introduction.}

\medskip
\noindent
The smooth dynamical system (vector field or diffeomorphism) is a Morse-Smale
dynamical system if:

\medskip
\noindent
1) It has finite number of critical elements (periodic trajectories for
diffeomorphisms, fixed points and closed orbits in the case of a vector field)
and all of them are non degenerated (hyperbolic);

\noindent
2) The stable and unstable manifolds of critical elements have transversal
intersections;

\noindent
3) Limit set of each trajectory is the critical element.

\bigskip

In papers [1] - [4] the topological classification of Morse-Smale dynamical
systems on 2-manifolds and in [4] - [7] on 3-manifolds are obtained.
Geteroclinic trajectories of the diffeomorphism are the trajectories laying in
the intersection of stable and unstable manifolds of critical elements of the
same indexes. Morse-Smale dynamic system is a gradient like if it hasn't contain
closed trajectories in the case of a vector field [8] and geteroclinic trajectories
in the case of diffeomorphism [9].

Two vector fields are topological equivalent, if there is a homeomorphism of
manifold to itself, which maps integral trajectories into integral
trajectories, keeping their orientations. This homeomorphism we shall call as
conjugated. Two diffeomorphism $f$, $g:M \rightarrow  M$ is topological
conjugated, if there is a homeomorphism $h: M \rightarrow  M$ such, that $hf =
gh.$

The purpose of this paper is to obtain topological classification of gradient
like Morse-Smale dynamical systems on 4-manifolds.

In section 2, we prove the criterion of topological equivalence of
vector fields on 4-manifolds. We construct invariant of field that
is 3-manifolds $N$ with
imbedded 2-spheres, framed links and surfaces. In section 3 and 4,
we use it for construction the diagram of vector fields and prove
the criteria of vector field equivalence using these diagrams. In
section 5, we construct the diagram of diffeomorphism. In section
6, we investigate when the diagram can be realized as diagram of
dynamical system.

\bigskip
\noindent
{\bf 2. Criterion of topological equivalence of vector fields.}

\medskip
Let $M$ be a closed 4-manifold, $X$ and $X'$ are Morse-Smale vector fields.

Let $a_1,...$, $a_k$ are fixed points of an index 0 of the field $X$, and
$a_1',...$, $a_k'$ of the field $X'$; $b_1...$, $b_n$ and $b_1',...$, $b_n'$
are fixed points of an index 1. Let $K$ is a union of stable manifolds of the
fixed points of index 0 and 1. We consider such tubular neighborhood $U(K)$ of
this union, which does not contain other fixed points and such, that each
trajectory has no more than one point of intersection with it. By $N =
\partial U(K)$ we denote the boundary of this neighborhood for the field $X$
and by $N'$ for the field $X'$. Then these boundaries are 3-manifold.

We denote the stable and unstable manifolds of the fixed point $X$ by $v(x)$
and $u(x)$. Let $S_i$ are surfaces, which are intersections of unstable
manifolds of the fixed points of an index 1 with manifold N. Then $S_i$ is a
set of not crossed spheres on manifold $N$, the supplement to which in $N$ is
three-dimensional spheres with cut out disks.

If $c_1,..,c_m$ are fixed points of an index 2, then the intersections
$v_i=v(c_i)\cap N$ form a set of the closed curve on manifold $N$. Similarly for
a field $X'$ on manifold $N'$ there is a set of the imbedded spheres and closed
curve.

A handle of index $k (k$-handle) $h^k=D^k\times D^{n-k}$ attached to manifold
$M$ with boundary is the union  $M \cup _f h^k$ along the boundary of $M$
according to an embedding

$$f: S^{k-1}\times D^{n-k}\rightarrow \partial M.$$

\noindent
$D^k\times \{0\}$ and $\{0\}\times D^{n-k}$ are called by core and cocore of the
handle. $S^{k-1}\times \{0\}$ and $\{0\}\times S^{n-k-1}$ are called by attached
and belt spheres of the handle $h^k$. Handles decompositions is a sequence of
imbeddings

$$M_0 \subset  M_1 \subset ...\subset  M_N =M$$

\noindent
such that $M_0$ is a union of $n$-disks (0-handles), $M_{i+1}$ is obtained from
$M_i$ by gluing handle. Handles decompositions are isomorphic, if there exist
a homeomorphism between the manifolds which maps the handles on the handles,
the cores to the cores and the cocores to the cocores.

In another way the manifold $N = \partial U(K)$ together with the imbedded spheres and closed
curves can be obtained if to consider the handle decomposition that is
associated with the given vector field. This is such handle decomposition, in
which to each handle there corresponds equally one fixed point (and on the
contrary), and the core of each handle lay on the stable manifold of the
appropriate fixed point, and the cocore lay on unstable one. Then $U(K)$ can be
considered as the union of the handles of the index 0 and 1, the imbedded
spheres as the belt sphere of the handles of the index 1, and closed curve as
the attached spheres of the handles of the index 2. The gluing of each 2-handle
is defined by a solid torus imbedded (neighborhood of the closed curve) in
the manifold N. It can be given by a closed curve and a parallel on torus or
simply by an integer (in a case if manifold $N$ is three-dimensional
sphere). This parallel is called by a framing.

We consider the intersections of the stable manifolds of the fixed points of an
index 3 with manifold $N$. These intersections are surfaces spanned on the closed
curves, that is such surfaces with boundaries, which interiors are imbedded in
$N$, and the boundary component coincide with the closed curves (one closed
curves can be contained in several component of the boundary). Thus the closed
curve so much time coincides with boundary components of a surface, as many
trajectories are in the intersections of unstable and stable manifolds of the
appropriate fixed points of an index 2 and 3.

\medskip
{\bf Theorem 1.} {\it Vector field $X$ is the topological equivalent to field $X'$
if and only if
there is a homeomorphism of manifolds $f:N \rightarrow N'$, which maps the
spheres in the spheres, the closed curves in the closed curves, keeping their
framing, and the surfaces in the surfaces.}

\medskip
Proof. Necessity. Let $\varphi $ is a conjugated homeomorphism  between fields $X$ and
$X'$. We set homeomorphism  $f$ as follows. For each point $x\in N$ we shall
consider a trajectory $g(x)$ which contain it. We define an image

$$f(x)=\varphi (\gamma (x))\cap N'.$$

As the unstable manifolds of the fixed points of index 1
map onto unstable manifolds of the fixed points of an index 1 by homeomorphism
$\varphi $, then the spheres map onto the spheres by homeomorphism $f$. Similarly,
the closed curves map onto the closed curves, and the surfaces onto the surfaces.
It is follows from existence of the homeomorphism $\varphi $ that appropriate
framings of the closed curve are equal.

\medskip
Sufficiency. Let there is a homeomorphism  $f: N \rightarrow  N'$. We
consider disks $D_i$, which lay on unstable integrated manifolds $u(b_i)$,
contain points $b_i$ and are limited by the spheres $S_i$. Then we can continue
homeomorphisms from boundaries $S_i$ of these disks up to homeomorphisms of
disks that map integrated trajectories onto integrated trajectories (as each
integrated trajectory, which has the intersections with a disk $D_i$, except
the fixed points, must intersect the boundary of the disk). Then the disks
$D_i$ are decompose a neighborhood $U(K)$ on 4-disks $D^4_j$, each of which
contain one fixed point of index 0.

Let us continue homeomorphisms from boundaries of these disks to an
interior. For it we consider three-dimensional spheres $S^3_j$, which are the
boundaries of such tubular neighborhoods of the fixed points of index 0, that
each trajectory from an interior $D^4_j$ crosses sphere $S^3_j$ in a unique
point. On each arch of the trajectory, which cross $D^4_j$ we introduce a new
parameterization in such a way that a point on $S^3_j$ corresponds to value of
the parameter $t = 0$, the point on the boundary $D^4_j$ to $t = 1$ and the
fixed point of the index 0 to $t = -1$. We demand that on each of the
intervals (-1; 0) and (0; 1) the parameters $t$ is proportional to the length
of the arch in the Riemann metrics that is fixed for all trajectories.
The homeomorphism of the boundaries of the disks $D^4_j$ sets a correspondence
between the trajectories inside this disk. A required homeomorphism maps a
point of each integrated trajectory to a point of the
appropriate integrated trajectory with the same parameter $t.$

For each closed curve $v_i$ we shall consider its tubular neighborhood
$U(v_i)$. Let $W_i$ is such neighborhood of the appropriate fixed point $c_i$,
which does not contain other fixed points and such, that boundary $\partial W_i
= U(v_i)\cup  V_i$ is a union of two solid torus. Thus the trajectories that have
the intersection with $W_i\backslash c_i$ go in $W_i$ through solid torus
$U(v_i))$ and leave it through solid torus $V_i$.

We cut solid torus $U(v_i)$ from manifold $N$ and glue solid torus $V_i$ to
obtain torus $\partial U(v_i)=\partial V$. Such operation we call
by spherical surgery along $v_i$. Fulfilling the spherical surgeries along
all $v_i$, we denote received manifolds by L. Thus between the solid torus $U(v_i)$ and
$V_i$ without middle circles (closed curve $v_i$ and $w_i$ on manifolds $N$ and
$L)$ there is a homeomorphism, which maps a point $x \in  N\backslash v_i$ to
$\gamma (x)\in L.$ Then homeomorphism of solid torus $U(v_i)$ and $U (v'_i)$
induce a homeomorphism  between $V_i\backslash w_i$ and $V'_i\backslash w'_i$.
Using equality of framing of the closed curve $v_i$ and $v'_i$ we can extend
this homeomorphism up to homeomorphism of the solid torus $V_i$ and $V'_i$. We
actually have constructed a homeomorphism between manifolds $L$ and $L'$. Thus
of surfaces from $N$ map onto imbedded 2-sphere in $L$, and 2-spheres map
onto surfaces. If homeomorphism from $N$ to $N'$ maps surfaces onto surfaces then
homeomorphism from $L$ to $L'$ maps 2-sphere onto 2-spheres.

We extend the homeomorphism from boundaries $W_i$ to their interior, and from
manifold $L$ to the other part of manifold $M$, just as we did with extension
of homeomorphism from $N$ to $U(K)$. Constructed homeomorphism of manifold $M$
will be required.

\bigskip
\noindent
{\bf 3.  Diagram of a polar vector field.}

\medskip
Let vector field is polar (i.e. with one source and one sink). In this
case, $K$ is a bouquet of circles and manifold $N$ is a connected sum
$\#_nS^1\times S^2$. We cut manifold $N$ by 2-spheres. The obtained manifold is
a 3-sphere without disks. Thus the closed curves are divided to arches, and
surfaces to surfaces with these arches and arch in 2-spheres as its boundaries.

So constructed three-dimensional sphere, together with the imbedded pairs of
2-spheres, with given homeomorphisms of spheres of one pair, arches and closed
curves with framing and surfaces we call by the diagram of a vector field. Two
diagrams are equivalent if there is a homeomorphism of three-dimensional sphere
of one diagram on three-dimensional sphere of other diagram, which maps

\medskip
1) pairs of 2-spheres in pairs of 2-spheres and commute with given homeomorphisms
of these spheres,

2) arch and closed curves in arches and closed curves and keep framing,

3) surfaces in surfaces.

\medskip
\noindent
If instead of surfaces we give the number of the fixed points (or handles) of
index 3, we obtain the Kirby diagram of manifold $M [10].$

\medskip
{\bf Proposition 1.} {\it Two polar gradient like Morse-Smale vector fields are
topological equivalent if and only if their diagrams are
equivalent.}

We will prove proposition 1 in general case.

\bigskip
\noindent
{\bf 4. Diagram of a vector field in a general case.}

\medskip
As above, we construct manifold $N$ with imbedded 2-dimensional spheres,
framing circles and surfaces with boundary on this circles. Cutting manifold
$N$ by spheres, we obtain 3-dimensional spheres with the imbedded pairs
of 2-dimensional spheres, as the boundaries of deleted disks, arches, which
connect them, and closed curves with the framing and surfaces with the
boundaries on this arch and closed curves.

The three-dimensional spheres, which is constructed in such way, together with
the pairs of embedded 2-dimensional spheres, with given homeomorphisms of
spheres of one pair, arches and closed curve with framing and surfaces with
boundaries on its we call by the diagram of the vector field. Two diagrams are
equivalent, if there is a homeomorphism  of three-dimensional spheres of one
diagram onto three-dimensional spheres of other diagram, which maps pairs of
2-dimensional spheres onto pairs of spheres, commuting with given homeomorphisms
of these spheres, arches and closed curves onto arches and closed curves, keeping
framing, and surfaces onto surfaces.

\medskip
{\bf Proposition 2.} {\it Two gradient like Morse-Smale vector fields are topological
equivalent if and only if their diagrams are equivalent.}

\medskip
Proof. Necessity. Required homeomorphism of three-dimensional spheres can be
obtained in result of restriction of the homeomorphism $f: N \rightarrow  N'$
on the correspondent three-dimensional spheres with deleted disks and extension
of it on these disks.

\medskip
Sufficiency. If there exist a homeomorphism between three-dimen\-sional spheres,
we consider its restriction on these spheres with deleted three-dimensional
disks, which are bounded by pairs of 2-dimensional spheres. Using that this
homeomorphism commute with given homeomorphisms of 2-dimensional spheres, we
can glue it in homeomorphism  of three-dimensional manifolds $f:N
\rightarrow N'$. Using the theorem 1 we shall receive necessary conditions of
the proposition.

\bigskip
\noindent
{\bf 5. Topological conjugation of diffeomorphisms.}

\medskip
Let $f: M \rightarrow  M$ is a gradient like Morse-Smale diffeomorphism.
Similarly to item 2 and 3 we shall construct the diagram of this
diffeomorphism. Then action of diffeomorphism $f$ on integrated manifolds of
the periodic saddle points induce the map of such sets:

\medskip
\noindent
1)   Three-dimensional spheres,

\noindent
2)   Pairs of 2-dimensional spheres, which is embedded in them,

\noindent
3)   Arches and circles with framing,

\noindent
4)   Surfaces,

\noindent
5)   Spheres with holes that obtained by cutting the manifold $L$ on
2-dimensional spheres, which are the images of surfaces with spherical surgery.

These maps we call by internal.

\medskip
{\bf Theorem 2.} {\it Two gradient like Morse-Smale diffeomorphisms $f$ and $g$ are
toplogicaly cojugated if and only if there exists a isomorphism of their
diagrams, which sets equivalence between them and commute with internal
maps.}

\medskip
Proof. The necessity of conditions follows from construction. Let's show
sufficiency. As it was done in the theorem 1, we construct homeomorphism $h$
between manifolds, which maps stable manifolds onto stable and unstable
onto unstable. Thus

$$h (f (U)) = g (h (U)),$$

\noindent
where $U$ is a part of stable or
unstable manifold, on which it is decomposed by other manifolds. Similarly, how
it was done for 3-dimensional manifolds [6,7], this
homeomorphism can be corrected up to required conjugated homeomorphism.

\bigskip
\noindent
{\bf 6. Realization of dynamical system with given invariant.}

\medskip
In this section we describe when the three-dimensional spheres, together with
the pairs of imbedded spheres, with given homeomorphisms of spheres of one
pair, arches and closed curves with framing and surfaces with boundary on its
are the diagram of some vector field.

We start with necessary conditions on surfaces. Since after spherical surgery
along the curves and arches with framing surfaces should became spheres, then
they are homeomorphic to spheres with deleted disks. Their intersections with
torus, which are the boundary of the curve or arches with framing, consist of
circles that set framing. The diagrams that satisfy to these conditions, we
call admissible.

\medskip
1) We shall consider the diagrams, in which there are no pairs of imbedded
spheres and surfaces. Each such diagram is the Kirby diagram. If it is the
diagram of a vector field, then this field have one source and sink and haven't
the fixed points of an index 1 and 3. Then after spherical surgery along the
framing link given by this diagram, the three-dimensional sphere should turn
out. Therefore Kirby diagram will be the diagram of a vector field in only case
when, when it is the diagram of three-dimensional sphere. From [10] it is follows,
that it will be in only case when from this diagram it is possible to receive
the empty diagram by additions of the closed curve and additions or removals
from the diagram of trivial knot with  framing $+1$ or $-1.$ Fenn and Rourke proved,
that these two movements can be replaced by one --- blow-up and its reverse ---
blow-down [11].

Nevertheless, these criteria do not give an algorithm, which check if the Kirby
diagram is the diagram of three-dimensional sphere or not. Such methods of
recognition of three-dimensional sphere having the Heegaard diagram are in [12]
and [13]. Using Rolfsen idea Christian Mercat show how the Heegaard diagrams
can be constructed from the Kirby diagrams.

\medskip
2) We show, as from any admissible diagram it is possible to obtain the Kirby
diagram (without embedded 2-dimensional spheres). In the beginning we get rid
from superfluous 1 and 3 handles (that is pairs of embedded spheres and
surfaces), which can be reduced together with additional 0 and 4 handles. If
the diagram consists of several three-dimensional spheres, then using
connectedness of the manifold \`I we can find the pair of 2-dimensional spheres
(1-handle) laying in different three-dimensional spheres. If we take thee
connected sum of three-dimensional spheres along this pair, we obtain the
diagram, with one three-dimensional sphere fewer than starting diagram. We
repeat this procedure until there will be one three-dimensional sphere.

Similarly it is possible to get rid from superfluous surfaces, consistently
deleting them. Thus at each stage the surface will be superfluous, if there not
exists the simple closed curve, which has transversal intersections with given
surfaces and does not have with others.

We replace the stayed pairs of spheres and surfaces with 2-handles. Thus we
shall receive the diagram from 1).

\bigskip
\noindent
\centerline{\bf REFERENCES}
}

\medskip
\noindent
1. S.H.Aranson, V.Z. Grines. Topological classification of flows on closed
2-manifolds $//$ Uspehi mat. nauk-41, No.1, 1986.- P.149-169.

\noindent
2. A.V.Bolsinov, A.A.Oshemkov, V.V.Sharko. On Classification of Flows on
Manifolds. I $//$ Methods of Funct. An. And Topology, v.2, N.2, 1996,
P.131-146.

\noindent
3. M. Peixoto. On the classification of flows on two-manifolds. In: Dynamical
Systems, edited by M.Peixoto.- Academic Press.- 1973.- P.389-419.

\noindent
4. G.  Fleitas Classification of Gradient like flows of dimensions two and
three $//$ Bol. Soc. Brasil. Mat.-9, 1975.- P.155-183.

\noindent
5. Ya.L.Umanski. The circuit of three-dimensional Morse-Smale dynamical system
without the closed trajectories $//$\ Sb. USSR, 230, N6,1976,- P.1286-1289.

\noindent
6. V.Z. Grines, H.H.Kalay. About topological classification of gradient like
diffeomorphisms on prime three-dimensional manifolds $//$ Uspehi mat. nauk,
$v.49$, No.2, 1994,- P.149-150.

\noindent
7.A.O. Prishlyak. Morse-Smale vector fields with finite number of the singular
trajectories on three-dimensional manifolds $//$ Doklady Nat.Ac.Sc. of
Ukraine, No.6, 1998.- P.43-47.

\noindent
8. Smale. S. On gradient dynamical system $//$ Ann. Math. 74, 1961, P.199-206.

\noindent
9. D.V. Anosov. Smooth dynamic systems 1 $//$ Results of a science and
engineering. Modern Problems of math. V.1, 1985, --- P.151-242.

\noindent
10. Kirby R.C. The Topology of 4-Manifolds. --- Lecture Notes 1374,
Springer-Verlag, 1989, -108p.

\noindent
11. R.Fenn, C. Rorke. On Kirby's calculus of links $//$ Topology,
v.18, No.1, 1979.- P.1-16.

\noindent
12. Rubinshtein J.N. An algorithm to recognize the 3-sphere $//$ Proc. ICM
(Zurich, 1994), Basel, 1995, P.601-611.

\noindent
13. Matveev S.V. A recognitial algorithm for 3-dimensional shpere (after
A.Tom\-pson)$//$ Mat. Sb. 186, 1995.- P.69-84.

\end{document}